\documentclass[12pt]{amsart}
\usepackage{amsmath,amstext,amsfonts,amsthm,amssymb}
\usepackage{eucal}
\usepackage{graphicx}
\renewcommand{\subsubsection}[1]{\addtocounter{subsubsection}{1}
{\ \\[3pt]\bf \thesubsubsection. \  #1} }

\swapnumbers

{  \theoremstyle{definition}

}
\newcommand{\am}{\operatorname{am}}

\newcommand{\cn}{\operatorname{cn}}
\newcommand{\dn}{\operatorname{dn}}
\newcommand{\sn}{\operatorname{sn}}

\newcommand{\br}{\text{{\bf r}}}

\newcommand{\lra}{\longrightarrow}

\newcommand{\bea}{\begin{eqnarray*}}
\newcommand{\eea}{\end{eqnarray*}}
\newcommand{\bean}{\begin{eqnarray}}
\newcommand{\eean}{\end{eqnarray}}

\newcommand{\BR}{\mathbb{R}}
\newcommand{\BZ}{\mathbb{Z}}

\begin{document}


\centerline{PENTAGRAMMA MIRIFICUM AND ELLIPTIC FUNCTIONS} 

\bigskip\bigskip

\centerline{(Napier, Gauss, Poncelet, Jacobi, ...)}

\bigskip\bigskip

\centerline{Vadim Schechtman}

\vspace{2cm}


FAUST. {\it Das Pentagramma macht dir Pein? 

Ei sage mir, du Sohn der H\"ole, 

Wenn das dich bannt, wie kamst du denn herein?}





\vspace{1cm}

\centerline{\bf Introduction}

\bigskip\bigskip

In this note we give an exposition of fragments from Gauss [G] where he discovered, 
with the help of some work of Jacobi,  
a remarkable connection between Napier pentagons on the sphere and Poncelet pentagons on the plane. 

This gives rise to a parametrization of the variety of Napier pentagons using the division by $5$ of elliptic functions. 

As a corollary we will find the classical five-term relation for the dilogarithm in a somewhat 
exotic disguise, cf. 6.4. 

The author gave talks on  these subjects in Toulouse on January, 2011 and in Moscow on May, 2011.

\bigskip\bigskip

\centerline{\bf \S 1. Napier's rules}

\bigskip\bigskip

John Napier of Merchiston (1550 - 1617) was a Scottish mathematician and astrologer. His main work is {\it Mirifici Logarithmorum Canonis Descriptio} (1614). The reader will find in [C] some historical 
remarks on what follows. 

{\bf 1.1.} For a spherical triangle with sides $a, b, c$ and angles $\alpha, \beta, \gamma$, 
$$
\cos c = \cos a \cos b + \sin a \sin b \cos \gamma
$$ 

Consider a right-angled spherical triangle with sides $a, b, c$ and angles $\alpha, \beta, \gamma = \pi/2$. 
Let us call its {\it parts} five quantities
$$
\tau =  (a, b, \pi/2 - \alpha, \pi/2 - c, \pi/2 - \beta)
$$
We consider them as positioned on a circle, i.e. ordered in the cyclic order. If we turn them 
we get the parts of another right-angled triangle: 
$$
\tau' = n(\tau) = (a', b', \ldots) = 
(b, \pi/2 - \alpha, \pi/2 - c, \pi/2 - \beta, a)
$$
({\it Napierian rotation}). Obviously $n^5 = 1$. 

Repeat this once more:
$$
g(\tau) := n^2(\tau) = (a'', b'', \ldots) =  
(\pi/2 - \alpha, \pi/2 - c, \pi/2 - \beta, a, b)
$$
In other words, 
$$
(\alpha'',b'',c'') = (\beta, \pi/2 - c, \pi/2 - a),
$$ 
this triangle is obtained by a reflection of the first one in the vertex $\beta$; call it {\it Gaussian reflection}. 


{\bf 1.2. Napier rules}: 

I. sine(middle part) = product of tangents of adjacent parts:
$$
\sin b = \tan a\cot\alpha
$$ 

II. sine(middle part) = product of cosines of opposite parts:
$$
\sin a = \sin \alpha \sin c
$$

\bigskip\bigskip



\centerline{\bf \S 2. Pentagramma Gaussianum}

\bigskip\bigskip

{\bf 2.1.} Let us draw a spherical right-angled triangle which we denote 
$P_3Q_1P_4$ with the angles $\angle Q_1P_3P_4 = p_3, \angle P_3Q_1P_4 = \pi/2, 
\angle Q_1P_4P_3 = p_4$ and sides $P_3Q_1  = \pi/2 - p_5, Q_1P_3 = \pi/2 - p_2, P_3P_4 = p_1$. 

Let us use an abbreviated notation $p' = \pi/2 - p$. So, the Napierian parts of this triangle are 
$$
\tau_1 = (p_2',p_5',p_3',p_1',p_4')
$$
Applying the Gaussian reflexion we get the triangles
$$
\tau_i := g^{i-1}(\tau_1) = (p'_{i+1},p'_{i+4},p'_{i+2},p'_{i+5},p'_{i+3})
$$
where we treat the indices modulo $5$. We denote the $i$-th triangle $P_{i+2}Q_iP_{i+3}$, 
where $\angle P_{i+2}Q_iP_{i+3} = \pi/2$, $p_i = P_{i+2}P_{i+3}$. 

So we get a spherical pentagon $P_1P_2P_3P_4P_5$; its characteristic property is $P_iP_{i+2} = \pi/2$, 
cf. the picture on p. 481 of [G]. 

Set
$$
\alpha_i = \tan^2 p_i
$$
Gauss' notation:
$$
(\alpha,\beta,\gamma,\delta,\epsilon) = (\alpha_1,\ldots,\alpha_5)
$$
Napier rules give:
$$
\cos p_i = \sin p_{i+1}\cos p_{i-1}
$$
$$
1 = \cos p_i\tan p_{i+2} \tan p_{i+3}
$$
It follows:
$$
(\gamma\delta,\delta\epsilon,\epsilon\alpha,\alpha\beta,\beta\gamma) = 
(\sec^2 p_1, \sec^2 p_2, \sec^2 p_3, \sec^2 p_4, \sec^2 p_5)
$$
Relations:
$$
1 + \alpha = \gamma\delta,\ 1 + \beta = \delta\epsilon,\ \text{etc.}
\eqno{(2.2.1)}
$$
Out of two quantities one can build up the remaining three, for example: 
$$
\beta = \frac{1 + \alpha + \gamma}{\alpha\gamma},\ \delta = \frac{1 + \alpha}{\gamma},\ 
\epsilon = \frac{1 +  \gamma}{\alpha},
$$
etc. and also 
$$
\gamma = \frac{1 + \alpha}{\alpha\beta - 1},\ \delta = \alpha\beta - 1,\ 
\epsilon = \frac{1 + \beta}{\alpha\beta - 1},
$$
etc. One has
$$
3 + \alpha + \beta + \gamma + \delta + \epsilon = \alpha \beta \gamma \delta \epsilon = 
\sqrt{(1+\alpha) (1+\beta) (1+\gamma) (1+\delta) (1+\epsilon)}
$$

{\bf 2.2.} Gauss' favorite example:  
$$
(\alpha, \beta, \gamma, \delta, \epsilon) = (9, 2/3, 2, 5, 1/3),\ 
\alpha \beta \gamma \delta \epsilon = 20
$$

{\bf 2.3. Regular pentagram.} If $\alpha = \beta = \ldots$ then 
$$
\alpha + 1 = \alpha^2,
$$
whence 
$$
\alpha = \frac{1 + \sqrt 5}{2}
$$
$$
\alpha^5 = \frac{11 + 5\sqrt 5}{2} = 11.0901699...
$$
Also $\alpha^2 = \sec^2 i$, so 
$$
\cos p_i = \alpha^{-1} = \frac{\alpha}{1 + \alpha} =  \frac{1 + \sqrt 5}{2}
$$
Note that
$$
\alpha^{-1} = \frac{\alpha}{1 + \alpha} =  \frac{-1 + \sqrt 5}{2} = 0.618033988...
$$

We have 
$$
\sin 2\pi/5 = 2s\cdot c
$$
where we have denoted $s = \sin \pi/5, c = \cos \pi/5$. On the other hand, 
$$
\sin 2\pi/5 = \sin 3\pi/5 = s(3 - 4s^2) = s(4c^2 - 1),
$$
whence $2c = 4c^2 - 1$, so
$$
\cos \pi/5 = \frac{1 + \sqrt 5}{4} = \frac{\alpha}{2}
$$ 
It follows that
$$
2\cos p_i\cdot \cos\pi/5 = 1 
$$
We set
$$
c' :=  \cos 2\pi/5 = \frac{- 1 + \sqrt 5}{4}
$$
Then
$$
cc' = 1/4,\ c - c' = 1/2
$$

\bigskip

{\it A cone }

\bigskip

{\bf 2.5.} Let $M$ be the centrum of the sphere. 
The vertices $P_i$ all 
lie on a quadratic cone whith vertex $M$. Namely, we have 
$$
\angle P_{k-1}MP_{k+1} = \pi/2
$$
So $MP_3$ is orthogonal to $MP_1$ and to $MP_5$. 

We define the coordinates in such a way that $M$ lies ar the origin,  
$$
P_3 = (1, 0, 0),\ 
P_1 = (0, 1, 0)
$$
Then
$$
P_5 = (0, \cos p_3, \sin p_3)
$$
and
$$
P_4 = (\cos p_1, 0, \sin p_1)
$$
Finally, the ray $MP_2\perp P_4MP_5$, so  
$$
P_2 = (\cos p_5, \cos p_4, - \cos p_3\sin p_5)
$$
The coordinates of $P_k$ satisfy the equation 
$$
\cos p_2\cdot(\cos p_1\cdot z - \sin p_1\cdot x)(\cos p_3\cdot z - \sin p_3\cdot y) = xy,
$$ 
or 
$$
z^2 - \sqrt\alpha xz - \sqrt\gamma yz - \frac{1 + \alpha + \gamma}{\sqrt{\alpha\gamma}} xy := 
$$
$$
= z^2 + p xz + q yz + r xy = 0
\eqno{(C)}
$$
Note that 
$$
r = - \beta\sqrt{\alpha\gamma}
$$

{\it Reduction to principal axes.}

{\bf 2.6.} In general, given a quadratic form defined by a symmetric matrix $A$, $v = (x, y, z)^*\mapsto 
v^tAv$, if we want to find an orthogonal matrix $B$ such that for $v = Bv'$, $v' = (x', y', z')^*$, 
$$
v^*Av = v^{\prime *}B^*ABv' = v^{\prime *}A'v'
$$
such that $A' = \text{diag}(G,G',G'')$ then $A' = B^*AB = B^{-1}AB$ and the numbers $G, G', G''$ 
are the eigenvalues of $A$. 

The columns of $B$ are eigenvectors of $A$. 

{\bf 2.7.} In our case, for the quadrics (C), the matrix is
$$
A = \left(\begin{matrix} 0 & r/2 & p/2\\ 
                         r/2 & 0 & q/2\\
                         p/2 & q/2 & 1\end{matrix}\right)
$$


The characteristic polynomial is 
$$
\det (t\cdot I - A) = t^3 - t^2 - \frac{p^2 + q^2 + r^2}{4} - \frac{r(pq - r)}{4}, 
$$
so the characteristic equation takes the form
$$
t(2t - 1)^2 = \omega(t - 1)
\eqno{(E)}
$$
where
$$
\omega := \alpha\beta\gamma\delta\epsilon
$$
This is our main equation.  

{\bf 2.8.} Let us investigate the real roots of $(E)$. We suppose that $\omega > 0$. The straight line 
$$
\ell:\ u = \omega(t - 1)
\eqno{(2.8.1)}
$$
allways intersects the cubic parabola
$$
P:\ u = t(2t - 1)^2
\eqno{(2.8.2)}
$$
at one negative point $t = G < 0$. 

On the other hand, $\ell$ intersects $P$ at two points $t = G', G'' > 1$ 
iff $\omega$ is greater than some critical value $\omega_0$; if $\omega = \omega_0$ then $\ell$ touches 
$P$ at $t = G' = G''$, if $\omega < \omega_0$ then there are no points of intersection of $\ell$ with 
$P$ other than $t = G$. 

The critical value is $\omega_0 = \alpha_0^5$ where 
$$
\alpha_0 = \frac{1 + \sqrt{5}}{2},
$$
cf. 2.3. In that case 
$$
G = - \alpha_0,\ G' = G'' = \alpha_0^2/2 = - c_0G,
$$
$$
c_0 = \cos\pi/5
$$
Thus, if $\omega\geq \alpha_0^5$ then      
$(E)$ has one negative root, $G$, and two positive roots, $G'$, $G''$ which coincide if $\omega = \alpha_0^5$.  

One has:
$$
G G' G'' = - \omega/4
$$
$$
(G - 1)(G' - 1)(G'' - 1) = - 1/4
$$
$$
(2G - 1)(2G' - 1)(2G'' - 1) = - \omega
$$

{\bf 2.9. Example.} For $\alpha\beta\gamma\delta\epsilon = 20$, $G = - 2.197, G' = 1.069,\ G'' = 2.128$.

\bigskip\bigskip



\centerline{\bf \S 3. Gauss' coordinates}

\bigskip\bigskip

{\bf 3.1.} Consider an ellipse
$$
E: \frac{x^2}{a^2} + \frac{x^2}{b^2} = 1,
$$
$a > b$. Inscribe $E$ into the circle $C$ with the centrum $O$ and radius $a$.
For a point $P = (x,y)\in E$, let $P' = (x,y')\in C$. 
The {\it eccentric anomaly (anomalie excentrique)} of $P$ is the angle $\phi = \angle XOP'$ where $X = (a,0)$. Then 
$$
x = a\cos \phi,\ y = b\sin \phi
$$

{\bf 3.2.} Return to the Napier pentagon $P_1\ldots P_5$.  
Let us draw a plane tangent to the sphere at the point of intersection with the axe of the cone 
and take the central projection of our pentagon on this plane. Let $R_i$ be the projection of $P_i$ 
to this plane. The points $R_i,\ i = 1, \ldots, 5,$ 
lie on an ellipse with semi-axes $\sqrt{-G/G'}$ and  $\sqrt{-G/G''}$. 

We put the coordinate axes $x, y$ of the plane along the axes of the ellipse. Let $R_i$ have the coordinates 
$(x_i, y_i)$. Let $\phi_i$ be the eccentric anomaly of $R_i$. We have
$$
x_i = g'\cos\phi_i,\ y_i = g''\sin\phi_i
\eqno{(3.2.1)}
$$
where 
$$
g' = \sqrt{- G/G'},\ g'' = \sqrt{- G/G''}
$$
Let $M$ be the centrum of the sphere, with coordinates $(x,y,z) = (0, 0, 1)$. 

{\it Attention}: these coordinates differ from 2.5.

Set 
$\psi_i = \angle P_{i}MP_{i+1} = \angle R_{i}MR_{i+1}$; $\br_i = (x_i,y_i,1)$.  

Then 
$$
\cos a_i = \frac{(\br_i,\br_{i+1})}{|\br_i||\br_{i+1}|} = \frac{x_{i}x_{i+1} + y_{i}y_{i+1} + 1}
{\sqrt{(x_i^2 + y_i^2 + 1)(x_{i+1}^2 + y_{i+1}^2 + 1)}} 
$$
and 
$$
\alpha_i := \tan^2 \psi_i = \sec^2 \psi_i - 1 = 
$$
$$
= \frac{(x_i - x_{i+1})^2 + (y_i - y_{i+1})^2 + (x_iy_{i+1} - y_ix_{i+1})^2}
{(x_{i}x_{i+1} + y_{i}y_{i+1} + 1)^2} = \frac{|\br_i\times \br_{i+1}|^2}{(\br_i,\br_{i+1})^2}
\eqno{(3.2.2)}
$$
For the future use, note the quantities
$$
\beta_i:= \sin^2 \psi_i = \frac{\alpha_i}{\alpha_i+1} = \frac{|\br_i\times \br_{i+1}|^2}{|\br_i|^2|\br_{i+1}|^2}
\eqno{(3.2.3)}
$$


{\bf 3.3.} We have
$$
\angle R_{i-1}MR_{i+1} = \pi/2
$$
whence
$$
x_{i-1}x_{i+1} + y_{i-1}y_{i+1} + 1 = 0
\eqno{(3.3.1)_i}
$$
Solving $(3.3.1)_{i-1}$ and $(3.3.1)_{i+1}$ for $x_i, y_i$, we get
$$
x_i = \frac{y_{i+2} - y_{i-2}}{x_{i+2}y_{i-2} - y_{i+2}x_{i-2}},\ 
y_i = \frac{x_{i-2} - x_{i+2}}{x_{i+2}y_{i-2} - y_{i+2}x_{i-2}}
\eqno{(3.3.2)}
$$

{\bf 3.4.} Next,
$$
\frac{x_ix_{i+1}}{2G' - 1} + \frac{y_iy_{i+1}}{2G'' - 1} + \frac{1}{2G - 1} = 0
\eqno{(3.4.1)}
$$
This will be proven in 3.6. 
It follows:
$$
x_i = - \frac{2G'-1}{2G-1}\cdot\frac{y_{i+1} - y_{i-1}}{x_{i-1}y_{i+1} - x_{i+1}y_{i-1}},
$$
$$
y_i =  \frac{2G''-1}{2G-1}\cdot\frac{x_{i+1} - x_{i-1}}{x_{i-1}y_{i+1} - x_{i+1}y_{i-1}}
\eqno{(3.4.2)}
$$ 

{\bf 3.5. Theorem} (Gauss, April 20, 1843). (a)
$$
\frac{\sin((\phi_{i-2} + \phi_{i+2})/2)}{\cos((\phi_{i-2} - \phi_{i+2})/2)} = 
\frac{G}{G''}\sin\phi_i;\ 
\frac{\cos((\phi_{i-2} + \phi_{i+2})/2)}{\cos((\phi_{i-2} - \phi_{i+2})/2)} = 
\frac{G}{G'}\cos\phi_i
\eqno{(3.5.1)}
$$
(b)
$$
\frac{\sin((\phi_{i-1} + \phi_{i+1})/2)}{\cos((\phi_{i-1} - \phi_{i+1})/2)} = 
\sqrt{\frac{G(G-1)}{G''(G''-1)}}\sin\phi_i = \frac{G(2G-1)}{G''(2G''-1)}\sin\phi_i
$$
$$ 
\frac{\cos((\phi_{i-1} + \phi_{i+1})/2)}{\cos((\phi_{i-1} - \phi_{i+1})/2)} = 
\sqrt{\frac{G(G-1)}{G'(G'-1)}}\cos\phi_i = \frac{G(2G-1)}{G'(2G'-1)}\cos\phi_i
\eqno{(3.5.2)}
$$

{\bf 3.6. Proof}, cf. [F]. (a) follows from (3.3.2) and (3.2.1), taking into account the formulas
$$
\sin a - \sin b = 2\cos((a+b)/2)\sin((a-b)/2
$$
$$
\cos a - \cos b = - 2\sin((a+b)/2)\sin((a-b)/2
\eqno{(3.6.1)}
$$


{\it Proof of} (3.4.1). 

It follows from (3.5.1) (replacing $i+2$ by $i$) that
$$
G^{\prime 2}\cos^2((\phi_i + \phi_{i+1})/2) + G^{\prime \prime 2}\sin^2((\phi_i + \phi_{i+1})/2) = 
G^2\cos^2((\phi_i - \phi_{i+1})/2).
$$
Using
$$
\cos^2 a = \frac{1 + \cos 2a}{2},\ \sin^2 a = \frac{1 - \cos 2a}{2}
$$
we get
$$
G^2\cos(\phi_i - \phi_{i+1}) + (G^{\prime\prime 2} - G^{\prime 2})\cos(\phi_i + \phi_{i+1}) = 
G^{\prime 2} + G^{\prime\prime 2} - G^2
$$
From (3.2.1) we get
$$
(G^{ 2} - G^{\prime 2} + G^{\prime\prime 2})G^{\prime}x_ix_{i+1} + 
(G^{ 2} + G^{\prime 2} - G^{\prime\prime 2})G^{\prime\prime}y_iy_{i+1} + 
$$
$$
+ (- G^{ 2} + G^{\prime 2} + G^{\prime\prime 2})G = 0
\eqno{(3.6.2)}
$$
The numbers $G, G', G''$ are the roots of 
$$
t(2t - 1)^2 = \omega(t - 1)
\eqno{(3.6.3)}
$$
It follows that
$$
G^2 + G^{\prime 2} + G^{\prime\prime 2} = \frac{1 + \omega}{2}
$$
Whence
$$
(- G^{ 2} + G^{\prime 2} + G^{\prime\prime 2})G = \biggl(\frac{1 + \omega}{2} - 2G^2\biggr)G = 
\frac{\omega}{2}\cdot \frac{1}{2G - 1}
$$
Similarly we compute the coefficients at $x_ix_{i+1}$ and at $y_iy_{i+1}$ of (3.6.2) and arrive 
at (3.4.1).

{\it Proof of} (b). The second equalities follow from (3.6.3). The first ones follow from (3.4.2) 
in the same manner as one has deduced (a) from (3.3.2). $\square$

\bigskip\bigskip



\centerline{\bf \S 4. Poncelet's problem and division of elliptic functions (Jacobi)}

\bigskip\bigskip

In this Section we describe the work of Jacobi [J2]. 
For a modern treatment of it see [GH] and references therein.  

{\bf 4.1.} One considers a circle with center $C$ of radius $R$, inside it a smaller circle with center $c$ of radius $r$, $a$ will be the distance between their centra. The line $cC$ intersects the bigger circle at a point $P$, so $|CP| = R$, $|cP| = R+a$.  

One takes a point $A_1$ on the bigger circle, draws the tangent to the smaller one till the intersection with the bigger one at the point $A_2$, and continues similarly. Denote 
$$
2\phi_i = \angle A_iCP
$$
Let us suppose for simplicity that $C$ is inside the smaller circle. Let $B$ be the point where 
$A_1A_2$ touches the smaller circle, $B'\in A_1A_2$ the base of the perpendicular from $C$, 
and $D\in cB$ the base ofthe perpendicular from $C$. Then
$$
DB = CB' = R\cos(\phi_2 - \phi_1)
$$
and 
$$
cD = a\cos(\phi_{2} + \phi_{1}),
$$
so
$$
r = cD + DB = R\cos(\phi_{2} - \phi_{1}) + a\cos(\phi_{2} + \phi_{1})
$$
It follows that  
$$
R\cos(\phi_{i+1} - \phi_{i}) + a\cos(\phi_{i+1} + \phi_{i}) = r
$$
or
$$
(R+a)\cos\phi_{i+1}\cos\phi_{i} + (R-a)\sin\phi_{i+1}\sin\phi_{i} = r
\eqno{(4.1.1)_i}
$$
Substracting $(4.1)_i - (4.1)_{i-1}$ and using
$$
\frac{\cos x - \cos y}{\sin y - \sin x} = \tan\biggl(\frac{x+y}{2}\biggr), 
$$
we get
$$
\tan((\phi_{i+1} + \phi_{i-1})/2)) = \frac{R - a}{R + a}\tan \phi_i
\eqno{(4.1.2)}
$$

{\bf 4.2. Jacobi elliptic functions.} Cf. [J1]. Fix a number $0\leq k < 1$ and define 
the period 
$$
K = \int_0^{\pi/2}\ \frac{dx}{\sqrt{1 - k^2\sin^2 x}} = 
\int_0^1\ \frac{dy}{\sqrt{(1 - y^2)(1 - k^2y^2)}}
\eqno{(4.2.1)} 
$$  
The function "amplitude" 
$$
\am:\ [0,K] \lra [0,\pi/2]
$$
is defined by $\phi = \am(u)$,  
$$
u = \int_0^\phi\ \frac{dx}{\sqrt{1 - k^2\sin^2 x}} = 
\int_0^{\sin \phi}\ \frac{dy}{\sqrt{(1 - y^2)(1 - k^2y^2)}}
$$
Thus 
$$
\am(K) = \pi/2
$$
We extend $\am$ to a function $\BR\lra \BR$ by $\am(-u) = - \am(u)$, $\am(u + 2K) = \am(u) + \pi$. 

We set
$$ 
\Delta\am u := \frac{d\am u}{du} = \sqrt{1 - k^2\sin^2\am u}
$$
Modern notations:
$$
\sn u = \sin\am u,\ \cn u = \cos\am u,\ \dn u = \Delta\am u
$$
Thus $\sn$ is an odd function and $\cn$ and $\dn$ are even.

{\bf 4.3. Addition formulas.}

$$
\sn(u+v) = \frac{\sn u\cn v\dn v + \cn u\sn v\dn u}{1 - k^2\sn^2 u\sn^2 v}
$$
$$
\cn(u+v) = \frac{\cn u\cn v - \sn u\sn v\dn u\dn v}{1 - k^2\sn^2 u\sn^2 v}
$$
$$
\dn(u+v) = \frac{\dn u\dn v - k^2\sn u\sn v\cn u\cn v}{1 - k^2\sn^2 u\sn^2 v},
\eqno{(4.3.1)}
$$
cf. [J1], \S 18. 

{\bf Main formula}:
$$
\cn(u-v) = \cn u\cn v + \sn u\sn v\dn (u-v)
\eqno{(4.3.2)}
$$

Next, 
$$
\sn(u+v) - \sn(u-v) = \frac{2\cn u\sn v\dn u}{1 - k^2\sn^2 u\sn^2 v}
$$
$$
\cn(u+v) - \cn(u-v) = - \frac{2\sn u\sn v\dn u\dn v}{1 - k^2\sn^2 u\sn^2 v},
$$
whence
$$
\tan((\am(x) + \am(y))/2) = \frac{\cn x - \cn y}{\sn y - \sn x} = 
$$
$$
\dn((x - y)/2)\tan\am((x + y)/2)
\eqno{(4.3.3)}
$$

{\bf 4.4.} If $\phi_n = \am(\phi + n t)$ then 
$$
\tan((\phi_0 + \phi_2)/2) = \Delta\am t\tan\phi_1
$$
So if 
$$
\Delta\am t = \frac{R-a}{R+a}
$$
then (4.1.2) is satisfied. 

We also have 
$$
\cos\phi_i\cos\phi_{i+1} + \sin\phi_i\sin\phi_{i+1}\cdot \sqrt{1 - k^2\sin^2\alpha} = 
\cos\alpha
$$
where $\alpha = \am t$. 

So we can find $k$ and $\alpha$ from the equations 
$$
\sqrt{1 - k^2\sin^2\alpha} = \frac{R-a}{R+a}
$$
and 
$$
\cos\alpha = \frac{r}{R+a}
$$
wherefrom
$$
k^2 = \frac{4Ra}{(R+a)^2 - r^2} = 1 - \frac{(R-a)^2 - r^2}{(R+a)^2 - r^2}
\eqno{(4.4.1)}
$$
It follows: 

{\bf 4.5. Theorem.} {\it If the process closes up after $n$ steps and $m$ full turns then, defining 
$k$  through} (4.4.1), {\it we have 
$$
\int_0^\alpha\ \frac{d\phi}{\sqrt{1 - k^2\sin^2\phi}} = \frac{m}{n} 
\int_0^\pi\ \frac{d\phi}{\sqrt{1 - k^2\sin^2\phi}}
$$}

\bigskip\bigskip



\centerline{\bf \S 5. Back to pentagramma}

\bigskip\bigskip

{\bf 5.1.} Let us return to (3.3.1). Denoting $\phi = \phi_i$, $\phi' = \phi_{i+2}$ and using (3.2.1) we get 
$$
g^{\prime 2}\cos\phi\cos\phi' + g^{\prime\prime 2}\sin\phi\sin\phi' = 1
$$
or
$$
\cos\phi\cos\phi' + \frac{G'}{G''}\sin\phi\sin\phi' = - \frac{G'}{G}
\eqno{(5.1.1)}
$$ 
Now compare this with the Main formula (4.3.2). We see that (5.1.1) will be satisfied if 
$\phi = \am(u), \phi' = \am(u + w)$, $\alpha = \am(w)$,  
$$
\dn w = \sqrt{1 - k^2\sin^2\alpha} = \frac{G'}{G''}
\eqno{(5.1.2)}
$$
and
$$
\cos\alpha = - \frac{G'}{G}.
$$
It follows that
$$
\cn w = - \frac{G'}{G}
\eqno{(5.1.3)}
$$
and
$$
k = \sqrt{\frac{G^{\prime - 2} - G^{\prime\prime -2}}{G^{\prime - 2} - G^{ -2}}}
\eqno{(5.1.4)}
$$
We can go backwards. 

{\bf 5.2. Theorem.} {\it For $0\leq k < 1$, let $K$ the corresponding complete elliptic integral} (4.2.1).   
{\it Define vectors in $ \BR^3$:  
$$
\br_j(k,u) = \biggl(\frac{\cn(u + 4jK/5)}{\sqrt{\cn(2K/5)}}, 
\frac{\sqrt{\dn(2K/5)}\sn(u + 4jK/5)}{\sqrt{\cn(2K/5)}}, 1\biggr), 
$$
$u\in \BR, j\in\BZ$. Then $\br_j(k,u) = \br_{j+5}(k,u)$. Set
$$
\alpha_j(k,u) = \frac{|\br_j(k,u)\times \br_{j+1}(k,u)|^2}{(\br_j(k,u),\br_{j+1}(k,u))^2}
$$
Then we have
$$
1 + \alpha_j(k,u) = \alpha_{j-2}(k,u)\alpha_{j+2}(k,u)
\eqno{(5.2.1)}
$$}



In the degenerate case  $k = 0$ we have
$$
|\br_j(0,u)|^2 = \sqrt 5,\ 
\alpha_j(0,u) = \frac{1 + \sqrt 5}{2}
$$
for all $j, u$. The ends of vectors $\br_j(0,u)$ form a regular plane pentagon, cf. 2.3.  

{\bf 5.3.} The idea of uniformising relations in the spherical (resp. hyperbolic) geometry by elliptic 
functions (cf. [S] for a nice review) plays an important role in Statistical Physics. Onsager uses it in his 1944 paper [O] on the solution 
of the Ising model. Baxter remarks that the same idea underlies the solution of his star-triangle 
relations, cf. [B], 7.13. 

Baxter cites a book [Gr] in this connection; apparently this is the very same 
"Greenhill's very odd and individual {\it Elliptic functions}" which, according to Littlewood cited 
by Hardy, was the Ramanujan's textbook, cf. [H], Ch. XII.

\bigskip\bigskip



\centerline{\bf \S 6. Dilogarithm five-term relation}

\bigskip\bigskip


{\bf 6.1.} {\it Euler dilogarithm}:
$$
Li_2(x) = \sum_{n=1}^\infty\ \frac{x^n}{n^2} = - \int_0^x\ \frac{\log(1-t)}{t}dt,
$$
$0\leq x \leq 1$. 

{\it Rogers dilogarithm}:  
$$
L(x) = Li_2(x) + \frac{1}{2}\log x\log(1-x) = -\frac{1}{2}\int_0^x\ 
\biggl(\frac{\log(1-t)}{t} - \frac{\log t}{1-t}\biggr) dt,
$$
$0 < x < 1$. 


{\bf 6.2. Theorem} (W.Spence, 1809). 
$$
L(x) + L(1-x) = \frac{\pi^2}{6}
\eqno{(6.2.1)}
$$
$$
L(x) + L(y) - L(xy) - L(x(1-y)/(1-xy)) - L(y(1-x))/(1-xy)) = 0, 
\eqno{(6.2.2)}
$$
$0 < x,y < 1$. 

{\bf 6.3.} One can reformulate (6.2.2) as follows (cf. [GT]): 
$$
L(x) + L(1-xy) + L(y) + L((1-y)/(1-xy)) + L((1-x)/(1-xy)) = \frac{\pi^2}{2}
$$
Define
$$
(b_1,\ldots,b_5) = (x,1-xy,y,(1-y)/(1-xy),(1-x)/(1-xy))
$$
and for an arbitrary $n\in \BZ$ define $b_n$ by periodicity $b_n = b_{n+5}$.  

Then
$$
b_{n-1}b_{n+1} = 1 - b_n,\ n\in \BZ.
$$
Define a new sequence $a_n$ by
$$
a_n = \frac{b_n}{1 - b_n},\ b_n = \frac{a_n}{1 + a_n}
$$ 
Then 
$$
a_{n-2}a_{n+2} = 1 + a_n
$$
Combining this with  5.2 and (3.2.3) we get

{\bf 6.4. Corollary.} {\it Let $\br_j(k,u)$ be as in} 5.2. {\it Set
$$
\beta_j(k,u) = \frac{|\br_j(k,u)\times \br_{j+1}(k,u)|^2}{|\br_j(k,u)|^2|\br_{j+1}(k,u)|^2}
$$
Then
$$
\sum_{j=1}^5 L(\beta_j(k,u)) = \frac{\pi^2}{2}
\eqno{(6.4.1)}
$$}
In the regular case $k = 0$ (6.4.1) takes the form 
$$
L(\alpha^{-1}) = \frac{\pi^2}{10}
\eqno{(6.4.2)}
$$
(Landen), cf. [K]. Here as usually $\alpha = (1 + \sqrt 5)/2$.

\bigskip

\bigskip\bigskip

\centerline{\bf References}

\bigskip\bigskip

[B] R.J.Baxter, Exactly solved models in statistical mechanics, Academic Press, 1982. 


[C] H.S.M.Coxeter, Frieze patterns, {\it Acta Arithm.} XVIII (1971), 297 - 304.  


[G] C.F.Gauss, Pentagramma Mirificum, Werke, Bd. III, 481 - 490; Bd VIII, 106 - 111.  

[F] R.Fricke, Bemerkungen zu [G], ibid. Bd. VIII, 112 - 117.  

[GT] F.Gliozzi, R.Tateo, ADE functional dilogarithm identities and integrable models, 
hep-th/9411203. 

[Gr] A.G.Greenhill, Applications of elliptic functions, 1892 (Dover, 1959). 

[GH] P.Griffiths, J.Harris, On Cayleys explicit solution to Poncelet's porism. 

[H] G.H.Hardy, Ramanujan, Cambridge, 1940 (AMS Chelsea, 1991). 

[J1] C.G.J.Jacobi, Fundamenta nova theoriae functionum ellipticarum.  

[J2] C.G.J.Jacobi, \"Uber die Anwendung der elliptischen Transcendenten auf ein bekanntes Problem 
der Elementargeometrie, {\it Crelles J.} {\bf 3} (1828)

[K] A.Kirillov, Dilogarithm identities, hept-th/9408113. 



[N] John Napier, Mirifici Logarithmorum canonis descriptio, Lugdini, 1619. 

[O] L.Onsager, Crystal statistics. I. A two-dimensioal model with an order-disorder transition, 
{\it Phys. Rev.} {\bf 65} (1944), 117 - 149. 


[S] J.Snape, Applications of elliptic functions in classical and algebraic geometry, Dissertation, Durham.  


\bigskip\bigskip

Institut de Math\'ematiques, Univ\'ersit\'e Paul Sabatier, 118 route de Narbonne, 31062 Toulouse, France

\end{document}